\documentclass{article}
\usepackage[english]{babel}
\usepackage{amsmath,amssymb,xcolor,tabularx}
\usepackage[utf8]{inputenc}
\usepackage[linesnumbered,ruled,vlined]{algorithm2e}
\usepackage{float, comment}
\usepackage[margin=1in]{geometry}
\usepackage{algorithmic}
\usepackage{xcolor}
\usepackage{caption}
\usepackage{booktabs}
\usepackage{array}
\usepackage{makecell}
\usepackage{abstract}
\usepackage[colorlinks, linkcolor=black, anchorcolor=blue, urlcolor=black, citecolor=black]{hyperref}
\usepackage{appendix}
\usepackage{bm}
\usepackage{authblk}
\usepackage{enumitem}
\usepackage[numbers,sort&compress]{natbib}
\usepackage{graphicx}
\usepackage{enumitem}
\setlist[itemize]{topsep=2pt, leftmargin=1.5em}
\usepackage{amsthm}
\theoremstyle{definition}
\newtheorem{remark}{Remark}

\newcommand{\assign}{:=}
\newcommand{\tmmathbf}[1]{\ensuremath{\boldsymbol{#1}}}
\newcommand{\tmop}[1]{\ensuremath{\operatorname{#1}}}

\newcommand{\tmtextbf}[1]{\text{{\bfseries{#1}}}}

\newcommand{\tmverbatim}[1]{\text{{\ttfamily{#1}}}}
\newenvironment{itemizedot}{\begin{itemize} }{\end{itemize}}

\newcounter{tmcounter}

\newcolumntype{Y}{>{\centering\arraybackslash}X}
\newcolumntype{L}[1]{>{\centering\arraybackslash}p{#1}}
\newcommand{\mailto}[1]{\href{mailto:#1}{#1}}

\begin{document}

\author{
  Haihao Lu \thanks{
  MIT, Sloan School of Management (\mailto{haihao@mit.edu});}
  \qquad
  Wanyu Zhang \thanks{
    Shanghai University of Finance and Economics 
    (\mailto{wanyuzhang@stu.sufe.edu.cn});}
}

\title{Enhanced PDHG for Linear Programming with Online Preconditioning}

\maketitle

\begin{abstract}
We present an online preconditioning technique for the primal-dual hybrid gradient (PDHG) algorithm for linear programming (LP). The method adaptively updates primal and dual preconditioners using an online optimization framework. To improve its practical performance, we introduce several algorithmic enhancements, including using normalized online loss functions and updating preconditioners infrequently. We implement the technique on top of vanilla PDHG and the GPU-based LP solver \tmverbatim{cuPDLP.jl} {\citep{lu2023cupdlp}}, and benchmark its performance on standard LP datasets. Our numerical experiments demonstrate that online preconditioning effectively reduces both iteration counts and overall solving time.
\end{abstract}

\section{Introduction}
There is growing interest in applying first-order methods (FOMs) to solve linear programming (LP) problems. One notable FOM-based LP solver is PDLP \cite{applegate2021practical, applegate2025pdlp}, which uses the primal-dual hybrid gradient (PDHG) algorithm. Unlike classical LP solvers that rely on matrix factorizations, PDLP operates primarily through matrix-vector multiplications, making it particularly well-suited for parallel computing and GPU acceleration. Its GPU implementations, \tmverbatim{cuPDLP.jl} \cite{lu2023cupdlp} and \tmverbatim{cuPDLP-C} \cite{lu2023cupdlp-c}, have attracted significant attention from both the optimization solver industry and major technology companies. These implementations represent a promising new direction in LP solving, and have demonstrated competitive, and in some cases superior, performance compared to commercial solvers in solving large-scale LPs.

To improve the practical performance of the PDHG algorithm, many numerical enhancements have been proposed. PDLP{\cite{applegate2021practical, applegate2025pdlp}} introduces several implementation enhancements, including adaptive restarts, preconditioning, adaptive stepsize, and feasibility polishing. \cite{lu2024restarted} incorporates Halpern iteration into PDLP and demonstrates accelerated computational performance. \cite{xiong2024role} develops the sublevel-set geometry condition numbers to analyze PDHG convergence and proposes a central-path Hessian-based rescaling method to improve these geometric measures and accelerate convergence. The success of PDHG in LP solving is rooted in these enhancements and motivates us to develop more effective and explainable algorithmic improvements.

In this work, we present an online preconditioning technique for PDHG in solving LPs and examine its practical performance. The technique is motivated by the online scaled gradient method
{\cite{gao2024gradient, chu2025provable, gao2025gradientmethodsonlinescaling}}
for smooth convex optimization, which formulates the problem of choosing preconditioners as an online decision-making problem (see Section~\ref{preliminary}). In our method, the primal and dual preconditioners of PDHG are updated by online gradient descent using feedback derived from the primal and dual updates, respectively. When specialized for LP solving, the method introduces minimal computational overhead—requiring only vector-vector operations—and is well-suited for parallelization on GPUs.

\begin{samepage}
The contributions of this work are summarized as follows:
\begin{itemizedot}
  \item We present an online preconditioning method for PDHG in LP solving. To improve its practical performance, we introduce several algorithmic enhancements, including integration with PDLP, the use of normalized online loss, and infrequent preconditioner updates (Section~\ref{online-pdhg}).
  
  \item We implement the online preconditioning technique on top of vanilla PDHG and the GPU-based LP solver \tmverbatim{cuPDLP.jl} {\cite{lu2023cupdlp}}, and empirically demonstrate its effectiveness on LP benchmark datasets (Section~\ref{experiments}).
\end{itemizedot}
\end{samepage}

\paragraph{Notations.} Let $\mathbb{S}^n_{+ +}$ denote the set of $n \times n$ symmetric positive-definite matrices. For $T \in \mathbb{S}^n_{+ +}$, we use $T^{1/2}$ to denote the unique symmetric matrix such that $T^{\frac{1}{2}} T^{\frac{1}{2}} = T$. We use $\|\cdot\|$ to denote the Euclidean norm for vectors and $\langle \cdot, \cdot \rangle$ for the corresponding inner product. The spectral norm of a matrix is denoted by $\|\cdot\|_2$. The notation $x \circ y$ denotes the element-wise (Hadamard) product of vectors. The projection operator onto a closed convex set $\mathcal{X}$ is written as $\tmop{proj}_{\mathcal{X}}(\cdot)$. For any vector $x$, $\tmop{Diag}(x)$ denotes the diagonal matrix with $x$ on its diagonal, and $I_n$ denotes the $n \times n$ identity matrix. The binary indicator function $\mathbb{I}\{\cdot\}$ equals $1$ if the condition holds and $0$ otherwise.

\section{Preliminaries}\label{preliminary}

\paragraph{Preconditioned PDHG for LP.} Consider solving a standard form LP:
\begin{eqnarray}
  \min_{x \in \mathbb{R}^n} & c^{\top} x &  \nonumber\\
  \text{s.t.} & A x = b &  \label{lp}\\
  & x \geq 0 &  \nonumber
\end{eqnarray}
where $A \in \mathbb{R}^{m \times n}$. The saddle-point formulation of problem \eqref{lp} is,
\begin{equation}
  \min_{x \in \mathbb{R}^n_+ } {\max_{\lambda \in \mathbb{R}^m}} L (x, \lambda) = c^{\top} x - \lambda^{\top} A x + b^{\top} \lambda . \label{pdhg-lp}
\end{equation}
To solve problem \eqref{pdhg-lp}, given primal and dual preconditioners $T \in \mathbb{S}^n_{+ +}, \Sigma \in \mathbb{S}^m_{+ +}$, the preconditioned PDHG iteration is:
\begin{align}
\begin{split}
  x^{k + 1} &= \operatorname{proj}_{\mathbb{R}^n_+} \left( x^k - T (c - A^{\top} \lambda^k) \right), \\
  \lambda^{k + 1} &= \lambda^k + \Sigma (b - A (2 x^{k + 1} - x^k)).
\end{split}
\label{ppdhg}
\end{align}
Various choices of $T$ and $\Sigma$ have been proposed in the literature:
\begin{itemizedot}
  \item Vanilla PDHG: $T = \tau \cdot I_n$, $\Sigma = \sigma \cdot I_m$, with the step sizes satisfying the condition $\tau \sigma \|A\|_{2}^{2}<1, \tau, \sigma > 0$.
  
  \item Pock-Chambolle preconditioning {\cite{pock2011diagonal}}: \( T = \operatorname{Diag}(\boldsymbol{\tau}) \), \( \Sigma = \operatorname{Diag}(\boldsymbol{\sigma}) \), where \( \boldsymbol{\tau} = (\tau_1, \ldots, \tau_n) \) and \( \boldsymbol{\sigma} = (\sigma_1, \ldots, \sigma_m) \). Given a parameter \( \beta \in [0, 2] \), the diagonal entries are defined as: 
   \[
  \tau_j = \frac{1}{\sum_{i = 1}^m |A_{i j}|^{2 - \beta}}, \quad \forall\, 1 \leq j \leq n; \quad \quad
  \sigma_i = \frac{1}{\sum_{j = 1}^n |A_{i j}|^{\beta}}, \quad \forall\, 1 \leq i \leq m.
  \]
  
  \item Ruiz preconditioning {\cite{ruiz2001scaling}}: \( T = \operatorname{Diag}(\boldsymbol{\tau}) \), \( \Sigma = \operatorname{Diag}(\boldsymbol{\sigma}) \), initialized as \( T = I_n \), \( \Sigma = I_m \). The method iteratively rescales the matrix \( A \) using the infinity norms of its rows and columns. Let \( \tilde{A} \) denote the scaled matrix, initialized as \( \tilde{A} = A \), then perform the following updates for \( K \) steps:
  \[
  \tilde A \leftarrow \Sigma^{\frac{1}{2}} \tilde A T^{\frac{1}{2}}; \quad
  \tau_j \leftarrow \frac{\tau_j}{\max_{1 \leq i \leq m} |\tilde A_{i j}|}, \quad \forall\, 1 \leq j \leq n; \quad
  \sigma_i \leftarrow \frac{\sigma_i}{\max_{1 \leq j \leq n} |\tilde A_{i j}|}, \quad \forall\, 1 \leq i \leq m.
  \]
Then, the preconditioners are scaled as $T = t \cdot \tmop{Diag} (\tmmathbf{\tau}), \Sigma = s \cdot \tmop{Diag} (\tmmathbf{\sigma})$, with scalars \( t, s > 0 \) chosen such that \( ts \|\tilde A\|_2^2 < 1 \) to ensure convergence.
\end{itemizedot}

\paragraph{Online preconditioning.}Consider solving an unconstrained minimization problem $\min_{x \in \mathbb{R}^n} f(x)$ using the preconditioned gradient descent method. At each iteration $k$, given a preconditioner $P_k \in \mathbb{R}^{n \times n}$, the update rule is:
\begin{align*}
  x^{k+1} &= x^k - P_k \nabla f(x^k).
\end{align*}
The online scaled gradient method (OSGM)~\cite{gao2024gradient, chu2025provable, gao2025gradientmethodsonlinescaling} formulates the problem of choosing $P_k$ as an online decision-making problem. Given the online feedback function $\ell_k(P_k)$, OSGM updates to $P_{k+1}$ by online gradient descent.
There are several choices of the online learning loss $\ell_k(P_k)$ in the literature:
\begin{itemizedot}
  \item Function value feedback: $\ell_k (P_k) = f (x^k - P_k \nabla f (x^k))
  - f (x^k)$.
  \item Hypergradient feedback {\cite{gao2024gradient}}: $\ell_k (P_k) =
  \frac{f (x^k - P_k \nabla f (x^k)) - f (x^k)}{\| \nabla f (x^k) \|^2}$.
\end{itemizedot}

Under function value feedback, the online gradient descent update becomes
\begin{equation}\label{eq:hdm}
P_{k + 1} = \operatorname{proj}_{\mathcal{P}} \left[ P_k - \alpha \nabla \ell_k (P_k) \right] = \operatorname{proj}_{\mathcal{P}} \left[ P_k + \alpha \nabla f (x^k - P_k \nabla f (x^k)) \nabla f (x^k)^{\top} \right],
\end{equation}
where $\mathcal{P}$ denotes the set of candidate preconditioners. Equation \eqref{eq:hdm} reduces to the hypergradient descent method
{\cite{baydin2018hypergradient}} when $\mathcal{P} \assign \{p \cdot I_{n} \mid p \in \mathbb{R}_{+}\}$. Hypergradient feedback is proposed by \cite{gao2024gradient} to incorporate normalization term $\|\nabla f(x^k)\|^2$. OSGM with hypergradient feedback has theoretical convergence guarantees when $f$ is convex and smooth.

\section{PDHG for LP with online preconditioning}\label{online-pdhg}

\subsection{Algorithm Derivation}\label{sec:basic-algo}

We extend the preconditioned PDHG update~\eqref{ppdhg} by allowing the primal and dual preconditioners to vary across iterations. Then we formulate the selection of preconditioners as an online learning problem.

At iteration \( k \), given the current iterate \( (x^k, \lambda^k) \), the algorithm selects preconditioners \( T_k \in \mathcal{P} \) and \( \Sigma_k \in \mathcal{D} \), where \( \mathcal{P} \) and \( \mathcal{D} \) denote the sets of candidate primal and dual preconditioners, respectively. A one-step preconditioned PDHG update is then performed:
\[
x^{k+1}(T_k) = \operatorname{proj}_{\mathbb{R}_+^n} \left( x^k - T_k (c - A^\top \lambda^k) \right), \quad
\lambda^{k+1}(\Sigma_k) = \lambda^k + \Sigma_k \left( b - A (2x^{k+1}(T_k) - x^k) \right).
\]
After the update, the algorithm receives feedback losses defined as:
\begin{align}
\begin{split}
  \ell^p_k(T_k) &= L(x^{k+1}(T_k), \lambda^k) - L(x^k, \lambda^k), \\
  \ell^d_k(\Sigma_k) &= L(x^k, \lambda^k) - L(x^k, \lambda^{k+1}(\Sigma_k)),
\end{split}
\label{def-loss}
\end{align}
where \( L(x, \lambda) \) denotes the Lagrangian defined in~\eqref{pdhg-lp}. The losses \( \ell^p_k(T_k), \ell^d_k(\Sigma_k) \) measure the change in the Lagrangian value due to the respective primal and dual updates.

To update the preconditioners, we apply online gradient descent (OGD) with step size \( \alpha > 0 \):
\[
T_{k+1} = \operatorname{proj}_{\mathcal{P}} \left[ T_k - \alpha \nabla \ell^p_k(T_k) \right], \quad
\Sigma_{k+1} = \operatorname{proj}_{\mathcal{D}} \left[ \Sigma_k - \alpha \nabla \ell^d_k(\Sigma_k) \right].
\]
Using the chain rule, the gradients of the loss functions can be written as:
\begin{align*}
  \nabla \ell^p_k(T_k) 
  &= - \left[ (c - A^{\top} \lambda^k) \circ \mathbb{I}\left( x^{k + 1/2}(T_k) \geq 0 \right) \right] (c - A^{\top} \lambda^k)^{\top}, \\
  \nabla \ell^d_k(\Sigma_k) 
  &= - (b - A x^k)(b - A (2 x^{k + 1} - x^k))^{\top},
\end{align*}
where \( x^{k + 1 / 2}(T_k) := x^k - T_k (c - A^\top \lambda^k) \). For computational efficiency, we restrict \( \mathcal{P} \) and \( \mathcal{D} \) to sets of nonnegative diagonal matrices throughout the remainder of this paper. Under this constraint, the preconditioners are updated in a coordinate-wise manner:
\begin{align*}
  T_{k+1} 
  &= T_k + \alpha \cdot \operatorname{Diag} \left( \left[ (c - A^\top \lambda^k) \circ \mathbb{I}\left( x^{k + 1/2}(T_k) \geq 0 \right) \right] \circ (c - A^\top \lambda^k) \right), \\
  \Sigma_{k+1} 
  &= \Sigma_k + \alpha \cdot \operatorname{Diag} \left( (b - A x^k) \circ (b - A (2 x^{k + 1} - x^k)) \right).
\end{align*}

\subsection{Practical algorithmic enhancements}\label{sec:pdlp-op}

In this section, we present several practical considerations for implementing the online preconditioning technique, including integration with PDLP, normalized online loss, and infrequent updates. The resulting one-step iteration incorporating these enhancements is summarized in Algorithm~\ref{alg:pdlp-op}.

\begin{algorithm}[htbp]
\caption{One-step iteration of PDHG with online preconditioning}
\label{alg:pdlp-op}
\begin{algorithmic}[1]
\REQUIRE Current iterate $(x^{k}, \lambda^{k})$, current preconditioners $T_{k}, \Sigma_{k}$, online learning rate $\alpha$, update frequency $\phi$, primal weight $\omega$, adaptive step size $\eta^{k}$
\ENSURE Next iterate $(x^{k+1}, \lambda^{k+1})$, updated preconditioners $T_{k+1}, \Sigma_{k+1}$

\STATE $(x^{k+1}, \lambda^{k+1}) \gets$ One step PDHG iteration with adaptive stepsize and primal weight \eqref{pdlp}
\IF{$k \bmod \phi \neq 0$}
    \STATE $(T_{k+1}, \Sigma_{k+1}) \gets (T_{k}, \Sigma_{k})$
\ELSE
    \STATE $(\nabla \ell^p_k(T_{k}), \nabla \ell^d_k(\Sigma_{k})) \gets$ Get online gradients via \eqref{pdlp-ongrad}
    \IF{using normalized online loss}
        \STATE $\nabla \ell^p_k(T_{k}) \gets \nabla \ell^p_k(T_{k}) / \|c - A^\top \lambda^{k}\|^2,\quad$ $\nabla \ell^d_k(\Sigma_{k}) \gets \nabla \ell^d_k(\Sigma_{k}) / \|b - A x^{k}\|^2$
    \ENDIF
    \STATE $T_{k+1} = \operatorname{proj}_{\mathcal{P}} \left[ T_k - \alpha \nabla \ell^p_k(T_k) \right], \quad$ $\Sigma_{k+1} = \operatorname{proj}_{\mathcal{D}} \left[ \Sigma_k - \alpha \nabla \ell^d_k(\Sigma_k) \right] $
\ENDIF
\end{algorithmic}
\end{algorithm}

\paragraph{Integration with PDLP.}PDLP {\cite{applegate2021practical}} incorporates several implementation tricks, including adaptive restarts, preconditioning, and adaptive stepsize. We now describe how to integrate the online preconditioning method into the PDLP framework.

Assume that after preprocessing (including preconditioning and presolving), the resulting LP has data $(c, b, A)$. Let $\eta^{k}$ denote the adaptive step size at iteration $k$, and $\omega$ the primal weight. Given online-learned preconditioners $T_{k}$ and $\Sigma_{k}$ (initialized as $T_{0} = I_n$ and $\Sigma_{0} = I_m$), the PDHG iteration with adaptive stepsize and primal weight becomes:
\begin{align}
\begin{split}
  x^{k+1}(T_{k}) &= \operatorname{proj}_{\mathbb{R}_+^n} \left( x^{k} - \frac{\eta^{k}}{\omega} T_{k} (c - A^\top \lambda^{k}) \right), \\
  \lambda^{k+1}(\Sigma_{k}) &= \lambda^{k} + \eta^{k} \omega \Sigma_{k} (b - A (2x^{k+1} - x^{k})).
\end{split}
\label{pdlp}
\end{align}

With feedback functions $\ell^p_k$ and $\ell^d_k$ defined in \eqref{def-loss}, the corresponding online gradients are given by:
\begin{align}
\begin{split}
  \nabla \ell^p_k(T_{k}) &= - \frac{\eta^{k}}{\omega} 
  \left[ (c - A^\top \lambda^{k}) \circ \mathbb{I}\left(x^{k + 1/2}(T_{k}) \geq 0\right) \right] 
  (c - A^\top \lambda^{k})^\top, \\
  \nabla \ell^d_k(\Sigma_{k}) &= - \eta^{k} \omega 
  (b - A x^{k}) (b - A (2 x^{k+1} - x^{k}))^\top.
\end{split}
\label{pdlp-ongrad}
\end{align}
We then update the preconditioners $(T_{k+1}, \Sigma_{k+1})$ via OGD using the online gradients in \eqref{pdlp-ongrad}. The adaptive restart scheme follows the same approach as that of  PDLP. When a restart is triggered, we initialize the next period's preconditioners by the most recent values.

\paragraph{\tmtextbf{Normalized online loss}.}Inspired by
{\cite{gao2024gradient}}, we normalize the online feedback by the squared norms $\| c - A^{\top} \lambda^k \|^2$ and $\| b - A x^k \|^2$. The resulting normalized loss functions are defined as
\begin{align*}
  \bar{\ell}^p_k(T_k) &= \frac{L(x^{k+1}(T_k), \lambda^k) - L(x^k, \lambda^k)}{\| c - A^{\top} \lambda^k \|^2}, \quad
  \bar{\ell}^d_k(\Sigma_k) = \frac{L(x^{k+1}, \lambda^k) - L(x^{k+1}, \lambda^{k+1}(\Sigma_k))}{\| b - A x^k \|^2}.
\end{align*}
However, this normalization requires additional computation of the squared norms.

\paragraph{Infrequent update.}The one-step update of preconditioners incurs two vector-vector products (and, if normalization is used, two additional norm computations). To reduce computational overhead, we adopt infrequent updates: we update the primal and dual preconditioners only once every $\phi$ PDHG iterations.

\section{Experiments}\label{experiments}

In this section, we evaluate the effectiveness of the online preconditioning technique. In Section~\ref{sec:exp-setup}, we describe the experiment setup. In Section~\ref{sec:exp1}, we illustrate the acceleration effect of online preconditioning on top of vanilla PDHG. In Section~\ref{sec:exp2}, we integrate the technique into the GPU-based LP solver \tmverbatim{cuPDLP.jl}~\cite{lu2023cupdlp}, and benchmark its performance on standard LP datasets.

\subsection{Experiment setup}\label{sec:exp-setup}

\paragraph{Benchmark datasets.} 
All experiments are conducted on instances from the Netlib collection~\cite{gay1985electronic} and the MIPLIB 2017 collection~\cite{gleixner2021miplib}. For the MIPLIB 2017 collection, we follow the setup in~\cite{applegate2021practical} and use 383 MIP relaxation instances.

\paragraph{Hyperparameters.} 
We set the relative optimality tolerance to $10^{-4}$ for all experiments. In Section~\ref{sec:exp1}, we impose an iteration limit of $5 \times 10^4$, and in Section~\ref{sec:exp2}, a time limit of 600 seconds is used. For the online learning rate scheduler, we use AdaGrad~\cite{duchi2011adaptive}. The online learning rate is selected on a per-instance basis from ${10^{-1}, 10^{-2}, 10^{-3}}$ in Section~\ref{sec:exp1}, and from ${10^{-5}, 10^{-6}, 10^{-7}}$ in Section~\ref{sec:exp2}. For infrequent updates, online preconditioners are updated every $\phi=20$ PDHG iterations. For other hyperparameters of \tmverbatim{cuPDLP.jl}, we use the default values.

\paragraph{Metrics.}In the benchmark experiments, the following performance metrics are reported:\begin{itemizedot}
  \item \#Opt: Number of \tmverbatim{OPTIMAL} instances solved within the time and tolerance limits.
  
  \item \#Iter (SGM10 / Mean): Shifted geometric mean (shift 10) and arithmetic mean of iteration counts over instances where all compared methods reach \tmverbatim{OPTIMAL}.
  
  \item Time (GM / Mean): Geometric mean and arithmetic mean of solving time (in seconds), over instances where all methods are \tmverbatim{OPTIMAL}.
  
  \item \#Imp. / Wors. Cases: Number of instances with reduced or increased iteration count compared to the baseline algorithm.
\end{itemizedot}

\paragraph{Computing environment.}All experiments are run on a computing cluster equipped with 8 NVIDIA A100 GPUs (80 GB VRAM each), running Ubuntu 22.04 with CUDA 12.5.

\subsection{Comparison with vanilla PDHG}\label{sec:exp1}

To understand the behavior of online preconditioning, we implement the online preconditioning technique on \tmverbatim{SimplePDLP.jl}\footnote{\url{https://github.com/MIT-Lu-Lab/SimplePDLP.git}}, a lightweight implementation of vanilla PDHG for LP solving that supports preconditioning. The experiments in this section use the default configuration of \tmverbatim{SimplePDLP.jl}, including a 10-step Ruiz rescaling followed by $L_{2}$-norm rescaling in the preprocessing phase. For online preconditioning, we apply the normalized online loss, and disable infrequent updates.

\begin{figure}[htbp]
\centering
\begin{tabular}{ccc}
  \includegraphics[width=0.3\textwidth]{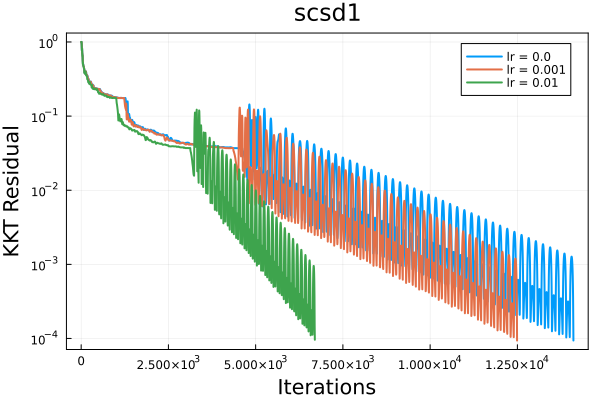} &
  \includegraphics[width=0.3\textwidth]{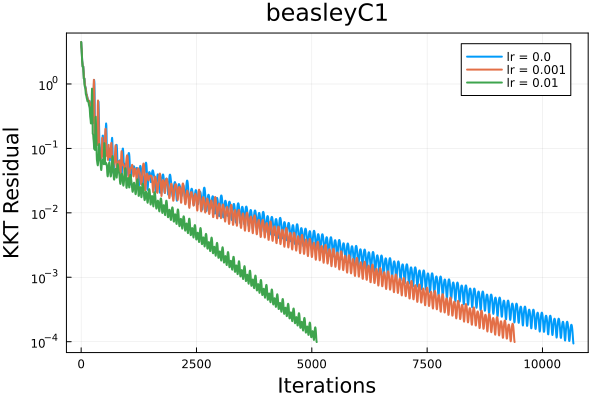} &
  \includegraphics[width=0.3\textwidth]{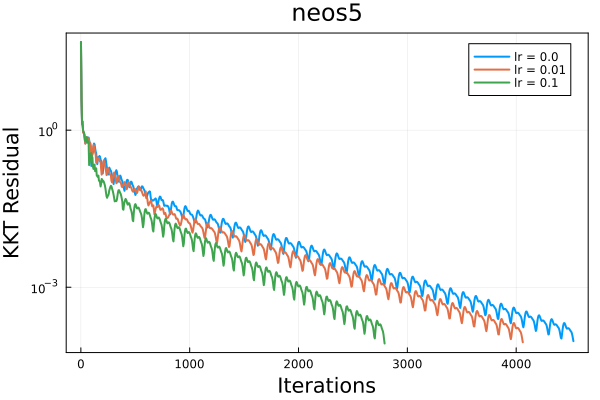} \\
  \includegraphics[width=0.3\textwidth]{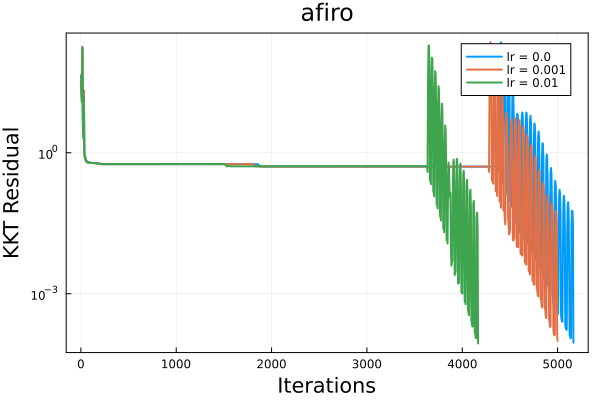} &
  \includegraphics[width=0.3\textwidth]{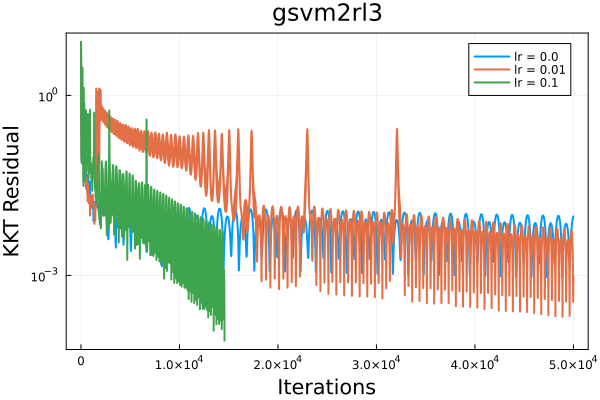} &
  \includegraphics[width=0.3\textwidth]{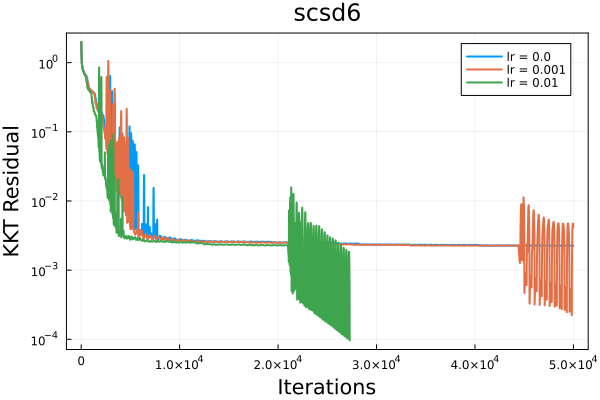} \\
  \includegraphics[width=0.3\textwidth]{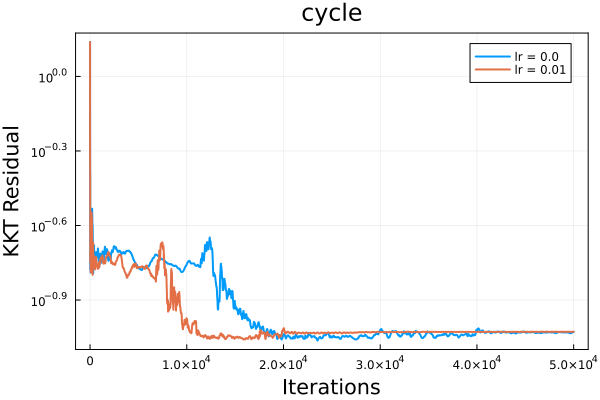} &
  \includegraphics[width=0.3\textwidth]{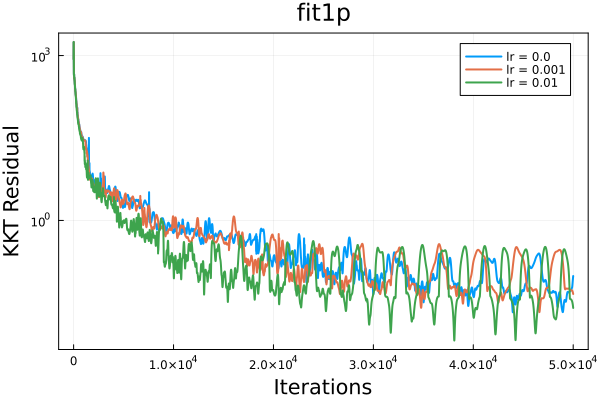} &
  \includegraphics[width=0.3\textwidth]{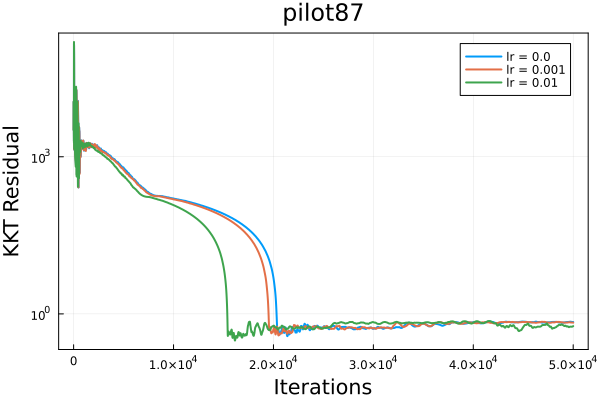}
\end{tabular}
\caption{The convergence behavior of vanilla PDHG with online preconditioning under different online learning rates (iteration limit $5 \times 10^4$, tolerance $10^{- 4}$)}
\label{fig:convergence}
\end{figure}

Figure \ref{fig:convergence} shows the convergence plots of PDHG with online preconditioning, under different online learning rates (lr), on instances from Netlib and MIPLIB dataset. When $\tmop{lr} = 0$ (the blue line), the method reduces to vanilla PDHG. The plots reveal the characteristic two-phase convergence behavior of PDHG~\cite{lu2024geometry}, featuring an initial relatively flat phase followed by a faster linear convergence phase. Based on these results, we categorize the observed convergence behaviors into three classes:

\begin{itemizedot}
  \item Vanilla PDHG is \tmverbatim{OPTIMAL} (\tmverbatim{scsd1}, \tmverbatim{beasleyC1},
  \tmverbatim{neos5}, \tmverbatim{afiro}): Online preconditioning improves the convergence rate, especially in the second phase.
  
  \item Vanilla PDHG reaches iteration limit, but becomes \tmverbatim{OPTIMAL} with online preconditioning (\tmverbatim{gsvm2rl3}, \tmverbatim{scsd6}): Online preconditioning enables a faster transition into the second convergence phase.
    
  \item Both vanilla PDHG and online preconditioning reach the iteration limit
  (\tmverbatim{cycle}, \tmverbatim{fit1p}, \tmverbatim{pilot87}).
\end{itemizedot}

\subsection{Comparison with \tmverbatim{cuPDLP.jl}}\label{sec:exp2}

We integrate the online preconditioning technique into \tmverbatim{cuPDLP.jl}, following the methodology described in Section~\ref{sec:pdlp-op}. We also investigate the effectiveness of two practical enhancements: normalization of online loss and infrequent updates of the preconditioner. These two binary options yield four variants of the online preconditioning method, summarized in the following table.

\begin{table}[htbp]
  \centering
  \begin{tabularx}{\textwidth}{|c|Y|Y|}
    \hline
     & Infrequent update - OFF & Infrequent update - ON \\
    \hline
    No normalization & \texttt{NoNorm-Freq} & \texttt{NoNorm-Infreq} \\
    \hline
    With normalization & \texttt{Norm-Freq}  & \texttt{Norm-Infreq}  \\
    \hline
  \end{tabularx}
\end{table}

The experiment results on the Netlib and MIPLIB datasets are summarized in Tables~\ref{netlib} and~\ref{miplib}, respectively. From these results, we make the following observations:
\begin{itemizedot}
  \item \tmtextbf{Online preconditioning reduces iteration counts.} All four variants of the online preconditioning technique yield improvements on iteration counts: the shifted geometric mean (SGM10) of iteration numbers decreases by approximately 10\% on both datasets, while the arithmetic mean decreases by around 30\% on MIPLIB.
  
  \item \tmtextbf{Normalization improves iteration efficiency.} Compared to their non-normalized counterparts, normalized variants further reduce SGM10 of iteration counts, regardless of update frequency. However, normalization introduces overhead due to the additional computation of vector norms, leading to increased solving time.
  
  \item \tmtextbf{Infrequent updates reduce runtime.} Using infrequent updates reduces solving time, with only a modest impact on iteration counts. On both datasets, variants with infrequent updates achieve lower geometric and arithmetic means of solving time.
  
  \item \tmtextbf{Benefits on large-scale instances.} The advantage of online preconditioning is more pronounced on large instances. We isolate instances for which all benchmark methods require more than 5000 iterations, and report results in Table~\ref{miplib-l}. On this subset, variant \tmverbatim{Norm-Infreq} achieves around 20\% reduction in solving time (GM).
\end{itemizedot}

\begin{remark}
The online learning rate is selected from ${10^{-5}, 10^{-6}, 10^{-7}}$ for each LP instance based on iteration counts. A tuning-free version of online preconditioning is under development.
\end{remark}

\begin{table}[h]
  \centering
  \caption{\label{netlib}Experiment results on Netlib dataset (time limit
  600 seconds, tolerance $10^{- 4}$)}
  \begin{tabularx}{\textwidth}{|L{2.5cm}|Y|Y|Y|Y|Y|Y|}
    \hline
    & \makecell{\#Iter\\(SGM10)} 
    & \makecell{\#Iter\\(Mean)} 
    & \makecell{Time\\(GM)} 
    & \makecell{Time\\(Mean)} 
    & \makecell{\#Imp.\\cases} 
    & \makecell{\#Wors.\\cases} \\
    \hline
    PDLP & 6356 & 15791 & 5.88 & 15.87 & -- & -- \\
    \hline
    \tmverbatim{NoNorm-Freq} & 5521 & \tmtextbf{13351} & \tmtextbf{4.86} & 12.08 & 66/95 & 20/95 \\
    \hline
    \tmverbatim{Norm-Freq} & 5482 & 13481 & 6.47 & 15.74 & \tmtextbf{68}/95 & 18/95 \\
    \hline
    \tmverbatim{NoNorm-Infreq} & 5380 & 13526 & 4.87 & \tmtextbf{12.05} & \tmtextbf{68}/95 & 18/95 \\
    \hline
    \tmverbatim{Norm-Infreq} & \tmtextbf{5266} & 13547 & 5.09 & 13.43 & \tmtextbf{68}/95 & \tmtextbf{16}/95 \\
    \hline
  \end{tabularx}
\end{table}

\begin{table}[h]
  \centering
  \caption{\label{miplib}Experiment results on MIPLIB dataset (time limit
  600 seconds, tolerance $10^{- 4}$)}
  \begin{tabularx}{\textwidth}{|L{2.5cm}|Y|Y|Y|Y|Y|Y|Y|}
    \hline
    & \makecell{\#Opt\\} 
    & \makecell{\#Iter\\(SGM10)} 
    & \makecell{\#Iter\\(Mean)} 
    & \makecell{Time\\(GM)} 
    & \makecell{Time\\(Mean)} 
    & \makecell{\#Imp.\\cases} 
    & \makecell{\#Wors.\\cases} \\
    \hline
    PDLP & 346 & 3693 & 25765 & 3.63 & 24.37 & -- & -- \\
    \hline
    \tmverbatim{NoNorm-Freq} & 341 & 3211 & 15476 & 4.82 & 21.92 & 183/338 & 69/338 \\
    \hline
    \tmverbatim{Norm-Freq} & 345 & \tmtextbf{3134} & \tmtextbf{15102} & 4.30 & 20.95 & \tmtextbf{201}/338 & \tmtextbf{50}/338 \\
    \hline
    \tmverbatim{NoNorm-Infreq} & 349 & 3214 & 18205 & \tmtextbf{3.07} & 17.33 & 187/338 & 51/338 \\
    \hline
    \tmverbatim{Norm-Infreq} & \tmtextbf{352} & 3187 & 15985 & 3.14 & \tmtextbf{15.74} & 185/338 & 62/338 \\
    \hline
  \end{tabularx}
\end{table}

\begin{table}[h]
  \centering
  \caption{\label{miplib-l}Experiment results on large instances (all
  methods haveiteration number $\geqslant$ 5000) of MIPLIB dataset (time limit
  600 seconds, tolerance $10^{- 4}$)}
  \begin{tabularx}{\textwidth}{|L{2.5cm}|Y|Y|Y|Y|Y|Y|}
    \hline
    & \makecell{\#Iter\\(SGM10)} 
    & \makecell{\#Iter\\(Mean)} 
    & \makecell{Time\\(GM)} 
    & \makecell{Time\\(Mean)} 
    & \makecell{\#Imp.\\cases} 
    & \makecell{\#Wors.\\cases} \\
    \hline
    PDLP & 26120 & 66103 & 24.43 & 61.91 & -- & -- \\
    \hline
    \tmverbatim{NoNorm-Freq} & 20503 & 38735 & 27.84 & 53.72 & 87/125 & 36/125 \\
    \hline
    \tmverbatim{Norm-Freq} & \tmtextbf{19322} & \tmtextbf{37718} & 26.49 & 52.24 & \tmtextbf{97}/125 & \tmtextbf{26}/125 \\
    \hline
    \tmverbatim{NoNorm-Infreq} & 20836 & 46178 & \tmtextbf{19.48} & 43.78 & 88/125 & 33/125 \\
    \hline
    \tmverbatim{Norm-Infreq} & 20327 & 40195 & 19.59 & \tmtextbf{39.38} & 90/125 & 32/125 \\
    \hline
  \end{tabularx}
\end{table}

\section{Conclusion}\label{conclusion}

In this work, we present an online preconditioning technique for the PDHG algorithm in LP solving. We implement this technique on top of vanilla PDHG and the GPU-based LP solver \tmverbatim{cuPDLP.jl}~\cite{lu2023cupdlp}, and incorporate normalization and infrequent updates to enhance practical performance. Our experiments demonstrate that online preconditioning can accelerate the convergence of PDHG and improve both iteration counts and overall solving time. Nonetheless, the current results rely on per-instance tuning of the online learning rate, and its theoretical convergence guarantees remain an open question. We leave both of these directions for future investigation.

\section*{Acknowledgments}

This technical report is part of the second author's undergraduate research conducted under the supervision of the first author. The authors would like to thank Wenzhi Gao for insightful discussions on the online scaled gradient method, and Dongdong Ge and Yinyu Ye for their support and encouragement throughout the project.

\bibliographystyle{plainnat}
\bibliography{ref}

\end{document}